\documentclass[10pt]{amsart}
\usepackage{amsmath, latexsym, 
}
\usepackage{amssymb,amsthm}
\usepackage{amscd,graphics}

\newtheorem{proposition}{Proposition}[section]
\newtheorem{lemma}[proposition]{Lemma}

\newtheorem{theorem}[proposition]{Theorem}

\theoremstyle{remark}

\newtheorem{remark}[proposition]{Remark}

\theoremstyle{definition}

\def\real{\mathbb{R}}

\def\complex{\mathbb{C}}
\def\supp{\mathrm{supp}}
\def\var{\mathrm{var}}
\def\sgn{\mathrm{sgn}}
\def\id{\mathrm{id}}

\def\JJ{\mathcal{J}}
\def\LL{\mathcal{L}}
\def\MM{\mathcal{M}}

\def \RR{\mathcal {R}}
\def\SS{\mathcal{S}}

\begin{document}
\title{On the susceptibility function of piecewise expanding interval maps}
\author{Viviane Baladi} 
\address{UMI 2924 CNRS-IMPA,  Estrada Dona Castorina 110,
22460-320 Rio de Janeiro, Brazil;
Permanent address: CNRS, UMR 7586, Institut de Math\'ematique de 
Jussieu, Paris}
\email{baladi@math.jussieu.fr}
\date{March 9 2007} 
\begin{abstract}
We study the susceptibility function 
$$\Psi(z)=\sum_{n=0}^\infty
\int z^n X (y) \rho_0(y) \frac {\partial} {\partial y} \varphi (f^n (y))\, dy$$
associated to the perturbation $f_t=f+tX\circ f$  of a piecewise expanding interval map $f$,
and to an observable $\varphi$.
$\Psi(1)$ is the formal
derivative (at $t=0$)
of the average   $\RR(t)=\int \varphi \rho_t\, dx$ of $\varphi$
with respect to the SRB measure of
$f_t$.
Our analysis is based on a spectral description of transfer
operators.
It gives in particular sufficient conditions on $f$, $X$, and $\varphi$ which guarantee
that $\Psi(z)$ is holomorphic in a disc of larger than one, or which ensure 
that a number may be associated to the divergent series $\Psi(1)$ .
We present
examples of $f$, $X$, and $\varphi$ 
so that $\RR(t)$ is not Lipschitz at $0$, and
we propose a new version of Ruelle's conjecture.
\end{abstract}
\thanks{Partially supported by ANR-05-JCJC-0107-01.
I  am grateful to  Dmitry Dolgopyat for important  remarks, and
to David Ruelle, who explained this problem to
me several years ago, and  shared  ideas on
his ongoing work on the nonuniformly hyperbolic case.
Artur Avila sketched the counterexample
of Theorem~\ref{counter}, provided Remark~\ref{AAA}, and made several useful comments.
I thank Gerhard Keller, who found mistakes in previous versions,
for helpful suggestions, and
Daniel Smania for very useful conversations who helped me formulate
Conjecture A and Remark~\ref{smania}}
\maketitle


\section{Introduction and main results}

Let us call SRB measure
for  a dynamical system $f: \MM \to \MM$, on a 
manifold $\MM$ endowed with Lebesgue measure, 
an $f$-invariant ergodic probability measures $\mu$ so that  the set
$
\{ x \in \MM\mid \lim_{n \to \infty}
\frac{1}{n} \sum_{k=0}^{n-1} \varphi(f^k(x))= \int \varphi \, d\mu \} 
$
has positive Lebesgue measure, for continuous
observables $\varphi$.
(Strictly speaking, this is the definition of a {\it physical measure,}  we refer
to \cite{Yo} for a discussion of the differences between physical and SRB measures.
For the purposes of this introduction, the distinction is not very important.)
If $f$ admits a unique SRB measure $\mu$, it is natural to ask
how $\mu$ varies when $f$ is changed. 
More precisely, one considers, for fixed $\varphi$, the
function
$
\RR(t)= \int \varphi \, d\mu_t
$, where $\mu_t$ is the SRB measure  (if it is well-defined) of $f_t=f+t X\circ f$.
Loosely speaking, we say that the SRB measure is differentiable
(or Lipschitz)  at
$f$ for $\varphi$ if $\RR(t)$ is differentiable at $0$.
(See \cite{Gal} for the relevance of this issue to nonequilibrium 
statistical mechanics.  Theorems 4 and 5 of \cite{Ke3}
show another setting where (Lipschitz) regularity of $\RR(t)$ is relevant.) 

If $f$ is a sufficiently
smooth uniformly hyperbolic diffeomorphism restricted to a transitive attractor,
Ruelle \cite{Ruok} (see also \cite{Ruok'}) proved that $\RR(t)$
is differentiable at $t=0$ and gave an explicit formula for
$\RR'(0)$.  Dolgopyat \cite{Dolgo} later showed that $\RR(t)$
was differentiable for a class of partially hyperbolic diffeomorphisms
$f$.  More recently, differentiability, together with a formula
for $\RR'(0)$, has been obtained for uniformly
hyperbolic continuous-time  systems (see \cite{BuL}
and references therein) and infinite-dimensional
hyperbolic systems (see \cite{JL} and references therein).

A much more difficult situation consists in studying
nonuniformly hyperbolic interval maps $f$, e.g. within the quadratic
family (not to mention higher-dimensional
dynamics such as H\'enon maps). For quadratic interval maps, one requires in
addition that the SRB measure be absolutely continuous with
respect to Lebesgue. It is well-known that
the SRB measure of $f_t$ may exist only for some
parameters $t$, although it is continuous 
in  a nontrivial subset of parameters (see \cite{RS}, \cite{Ts}). 
In this setting,  Ruelle (\cite{Rupr}, \cite{Ru3})
has outlined a program, replacing
differentiability by differentiability
in the sense of Whitney's extension theorem, and proposing
$\Psi(1)$ with
\begin{equation}\label{intro}
\Psi(z)=\sum_{n=0}^\infty
\int z^n X (y) \rho_0(y) \frac {\partial} {\partial y} \varphi (f^n (y))\, dy\, ,
\end{equation}
the ``susceptibility function,"
\footnote{Since $\Psi(e^{i\omega})$
is the Fourier transform of the ``linear response" 
\cite{Ru1ph}, it is natural to consider the variable
$\omega$, but we prefer to work with the variable $z=e^{i\omega}$.}
as a candidate 
for the derivative. Beware that $\Psi(1)$ needs to be suitably
interpreted: It could be simply the value at $1$ of a meromorphic extension
of $\Psi(z)$ such that $1$ is not a pole, but also a number associated to
the -- possibly divergent -- series obtained by setting $z=1$ in
(\ref{intro}), by some (yet undetermined) summability method.
Formal arguments (see \cite{Ru1ph} and Appendix ~\ref{susc}) justify the choice of $\Psi(1)$,
which Ruelle \cite{Rupr} calls ``the only reasonable formula one can write."
For several nonuniformly hyperbolic
interval maps  $f$ admitting a finite Markov partition
(i.e., the critical point is preperiodic),
although $\Psi(z)$ has a pole (or several poles)
inside the open unit disc, it extends meromorphically to
a disc of radius larger than $1$ and is holomorphic at $z=1$
(\cite{Ru3}, \cite{JR}). The relation between $\Psi(1)$
and (Whitney)
differentiability of $\RR(t)$ for such maps has not been established.
The case of nonrecurrent critical points is being
investigated \cite{Ruwip}.

\smallskip

Our goal here is much more modest: We consider unimodal interval
maps $f$ which are piecewise {\it uniformly} expanding, i.e., $|f'|> 1$ (except at
the critical point). In this case, existence of the SRB measure of
all perturbed maps $f_t$ is guaranteed, and it is  known 
that $\RR(t)$ has modulus of continuity
$|t|\ln |t|$ (we refer to the beginning of Section~\ref{two} for
more details and references). Our intention was to understand the analytic
properties of $\Psi(z)$ for perturbations $f+tX\circ f$ of such maps, and
to see if they could be related to the differentiability 
(or lack of differentiability) of $\RR(t)$.
Our results are as follows (the precise setting is described in
Section~\ref{two}):

We prove (Proposition~\ref{naive}) that $\Psi(z)$ is always
holomorphic in the open unit disc. When the critical point
is preperiodic of eventual period $n_1\ge 1$, we show that 
(Theorems~\ref{M1} and ~\ref{M2})  $\Psi(z)$ extends meromorphically to
a disc of radius larger than one, with possible poles
at the $n_1$th roots of unity, and we give sufficient conditions
for the residues of the poles to vanish. When the critical point is not
periodic, $\Psi(z)$ appears  to be rarely holomorphic at
$z=1$. Nevertheless, we have a candidate $\Psi_1$ for the value of the
possibly divergent series $\Psi(1)$, under 
the condition that the ``weighted total
jump" $\JJ(f,X)$ defined in (\ref{saltusid}) vanishes (Proposition~\ref{prop1}).
The tools for these results are  transfer operators $\LL_0$ and
$\LL_1$ introduced in Section~\ref{two} (these operators were also used by Ruelle 
\cite{Ru3}). A key ingredient is a decomposition (Proposition~\ref{decc})
of the invariant density of $f$ into a smooth component and a
``jump"  component (this 
was inspired by Ruelle's work \cite{Ruwip} in the nonuniformly
hyperbolic case).

Finally, we give 
examples of interval maps and observables for which $\RR(t)$ is not Lipschitz.
\footnote{After this paper was written, Carlangelo Liverani mentioned to us
that Marco Mazzolena \cite{MM} independently constructed examples of families
$f_t$ such that $\RR(t)$ is not Lipschitz.}
Applying Theorem~\ref{M1} to
these examples we get that $\Psi(z)$  has a pole at $z=1$. 
The ``weighted total jump" $\JJ(f,X)$ associated to these 
examples is nonzero.

\smallskip

In view of our results, we propose
to reformulate Ruelle's conjecture as follows:

\smallskip
{\bf Conjecture A.}
Let $f$ be either a mixing, piecewise expanding, piecewise
smooth unimodal interval map such that the critical point is
not periodic, or a mixing smooth Collet-Eckmann
unimodal  interval map with nonflat critical point.
Let $f_t=f+X_t \circ f$ be a smooth perturbation (with $X_0=0$) corresponding
to a smooth $X =\partial_t X_t|_{t=0}$
such that each $f_t$ is topologically conjugated to
$f$. Then $\RR(t)$ is differentiable at $0$ for
all smooth observables $\varphi$, and $\RR'(0)= \Psi(1)$ (the
infinite sum being suitably interpreted).

\smallskip

The above conjecture is interesting only if there are
examples satisfying the assumptions and for which the conjugacy
between $f$ and $f_t$ is not smooth. We  \cite{BS} expect this to
be true and that the condition $\JJ(f,X)=0$  is related to the existence
of  a topological conjugacy between $f$ and $f_t$ (see Remark ~\ref{smania}).

\medskip

For general perturbations of piecewise expanding maps, our counter-examples
show that the (previously known) property that
$\RR(t)$ has modulus of continuity $|t|\ln|t|$ cannot be improved.
For nonuniformly expanding maps, we propose:

{\bf Conjecture B.}
Let $f$ be  a mixing smooth Collet-Eckmann
unimodal  interval map, with nondegenerate critical point $c$
(i.e. $f''(c)\ne 0$). Then, for any smooth  $X$, 
and any $C^1$ observable $\varphi$, the function $\RR(t)$ is
$\eta$-H\"older at $0$, in the sense of Whitney
over those $t$ for which $f_t$ is Collet-Eckmann, for any $\eta < 1/2$. 
\smallskip

For critical points of order $p\ge 3$ we expect  
that the condition $\eta < 1/2$ should be replaced by $\eta < 1/p$.
We expect Conjectures A and B to be essentially optimal.

\bigskip
\section{Setting and spectral properties of the transfer operators} 
\label{two}

In this work, we consider a continuous 
$f:I\to I$ where $I=[a,b]$, with:
  \renewcommand{\labelenumi}{(\roman{enumi})}
\begin{enumerate}
\item
$f$ is strictly increasing
on $I_+=[a,c]$, strictly decreasing on $I_-=[c,b]$
($a<c<b$),
\item\label{last}
for $\sigma=\pm$, the map
$f|_{I_\sigma}$ extends to a  $C^3$ map
on a neighbourhood of $I_\sigma$, and
 $\inf |f'|_{I_\sigma}|  > 1$;
\item\label{nonper}
$c$ is  not periodic under $f$;
\item\label{topmix}
$f$ is topologically mixing on
$[f^2(c), f(c)]$.
\end{enumerate}
\smallskip

The point $c$ will be called the {\it critical point} of $f$.
We write $c_k=f^k(c)$ for $k \ge 0$.

For a function $X: \real \to \real$,  with
$\sup |X| \le 1$,
so that $X|_{f(I)}$
extends to a $C^2$ function in a neighbourhood
of $f(I)$ and $X'$ is of bounded variation \footnote{A prime denotes derivation, a priori in the sense of distributions.}  and supported in $[a,b]$,
we shall consider 
{\it the additive perturbation}
\footnote{Sometimes we only consider one-sided perturbations, i.e., $t \ge 0$ or $t\le 0$.}
\begin{equation}
f_t(x)=f(x)+t X(f(x))\, ,
\quad | t| \mbox{ small.}
\end{equation} 
More precisely, we take $\epsilon>0$  so that
(i) and (ii)  hold for all $f_t$ with $|t|<\epsilon$,
except that $f_t|_{I\sigma}$ may only extend to a
$C^2$ map. Then we {\it assume} that $f$ and $X$ are
such that, up to taking
perhaps smaller $\epsilon$, we have
$\sup_{|t|<\epsilon} f_t(c)\le b$ 
and $\inf_{|t|<\epsilon}\min( f_t(a), f_t(b))\ge a$, so
that each $f_t$ maps $I$ into itself.
Then each $f_t$ admits an absolutely
continuous invariant probability measure, with a density
$\rho_t$ which is of bounded variation \cite{LY}.
There is only one such measure \cite{LiY} and it is ergodic.
In fact, assumption (iv) implies that it is mixing.
(We refer to the introduction of \cite{KL} for an
account  of the use of bounded variation spaces,
in particular references to the work of Rychlik and Keller.
The bibliography there, together with that in Ruelle's book
\cite{Rub}, give a fairly complete picture.)
By construction, each $\rho_t$ is continuous
on the complement of the at most countable set
$C_t=\{ f^k_t(c)\, , \, k\ge 1\}$,
and it is supported in $[f^2_t(c), f_t(c)]\subset [a,b]$
(we extend it by zero on $\real$). 
In addition, denoting by $| \varphi|_1$  the $L^1(\real,\mbox{Lebesgue})$ norm of $\varphi$,
assumption (iii) implies by
\cite[Prop. 7]{KL} (see (\ref{LYu}) below and also \cite[Remark 5]{KL}) that  
\begin{equation}\label{claim'}
 |\rho_t - \rho_0|_1 = 0(|t| \ln |t|) \, .
\end{equation}

\bigskip

We next define the transfer operators $\LL_0$ and $\LL_1$,
with $\LL_1$ the ordinary Perron-Frobenius operator,
and show that $\LL_1$ is ``the derivative of $\LL_0$."
In order to make this precise we need more notation.
Recall that a point $x$ 
is called regular  for a function $\phi$ 
if $2\phi(x)=\lim_{y \uparrow x} \phi(y)+\lim_{y\downarrow x}\phi(y)$.
If $\phi_1$ 
and $\phi_2$ are (complex-valued) functions of
bounded variation on $\real$
having at most regular discontinuities, the Leibniz formula
says that $(\phi_1 \phi_2)'=\phi_1'\phi_2+\phi_1 \phi_2'$, where both
sides are a priori finite measures.
Define $J:=(-\infty , f(c)]$ and $\chi:\real \to \{0,1,1/2\}$ by
$$
\chi(x)=\begin{cases}
0& x \notin J\\
1& x \in \mbox{int}\,  J\\
\frac{1}{2}& x=f(c)\, .
\end{cases}
$$
The two inverse branches of  $f$,
a priori defined on $[f(a),f(c)]$ and $[f(b),f(c)]$,
may be extended to $C^3$ maps $\psi_{+}: J \to (-\infty,c]$
and $\psi_{-}: J \to [c,\infty)$, with $\sup |\psi_{\sigma}'|<1$ for $\sigma=\pm$.
(in fact there is a $C^3$ extension of $\psi_\pm$ in a small neighbourhood
of $J$.)
The map $\psi_{+}$ has a single fixed point
$a_0 \le a$.

\smallskip

We can now introduce two linear operators:
\begin{equation}\label{L0}
\LL_{0} \varphi(x):=
\chi(x) \varphi(\psi_{+}(x))-\chi(x)  \varphi(\psi_{-}(x))\, ,
\end{equation}
and
\begin{equation}\label{L1}
\LL_{1} \varphi(x):=
\chi(x)
\psi'_{+}(x)\varphi(\psi_{+}(x))
+\chi(x)
|\psi'_{-}(x)| \varphi(\psi_{-}(x))\, .
\end{equation}
Note that $\LL_{1}$
is  the usual (Perron-Frobenius) transfer operator for $f$,
in particular, $\LL_{1} \rho_0=\rho_0$ and $\LL_{1}^*(\mbox{Lebesgue}_\real)
=\mbox{Lebesgue}_\real$.
The operators $\LL_{0}$ and $\LL_{1}$  both act boundedly on 
the Banach space
$$
BV=BV^{(0)}
:=\{ \varphi : \real \to \complex \,,
\var (\varphi) < \infty\, , \supp (\varphi) \subset [a_0 , b]\} / \sim\, ,
$$ 
endowed with the norm $\| \varphi\|_{BV} =\inf _{\phi \sim \varphi} \var (\phi)$,
where $\var (\cdot)$ denotes the total
variation and $\varphi_1 \sim \varphi_2$ if the
bounded functions $\varphi_1$, $\varphi_2$ 
differ on an at most countable set. 

The following lemma indicates that
$BV$ is  the ``right space" for $\LL_{1}$,
but is not quite good enough  for $\LL_0$:

\begin{lemma}\label{spBV}
There is $\lambda <1$ so that 
the
essential spectral radius of $\LL_{1}$ on $BV$ is $\le \lambda$, while
$1$ is a maximal eigenvalue of $\LL_{1}$, which is simple, for the eigenvector
$\rho_0$. There are no other eigenvalues of $\LL_1$ of modulus $1$ on $BV$.

The essential spectral radius of 
$\LL_0$ on $BV$ coincides with its spectral radius,
and they are equal to $1$. 
\end{lemma}

\begin{proof}
For the claims on $\LL_1$, we refer  e.g. to \cite[\S 3.1--3.2]{Bb}
and references therein to works of Hofbauer, Keller, Baladi, Ruelle
(see also Appendix~ \ref{uLY}). In fact, we may take any
$\lambda \in (\sup_{x} (|f' (x)|^ {-1}), 1)$.

The
essential spectral radius of $\LL_{0}$ on $BV$
is equal to $1$ (see e.g. \cite[\S 3.2]{Bb}, and in particular the result
of \cite{Ke2} for the lower bound).
It remains to show that there are no
eigenvalues of modulus larger than $1$.
Now, $z$ is an eigenvalue of modulus $>1$
of $\LL_{0}$ on $BV$ if and only if
(see e.g.  \cite{Rub}) $w=1/z$ is a pole
of
$$
\zeta(w)=
\exp \sum_{n \ge 1} \frac{w^n}{n}
\sum_{f^n(x)=x} \sgn (f^n)'(x) \, .
$$
However, since $f$ is continuous, we have that $|\sum_{f^n(x)=x} \sgn (f^n)'(x) |\le 2$
for each $n$, so that $\zeta(w)$ has no poles in the open unit disc.
\end{proof}

To get finer information on $\LL_0$, we consider the smaller Banach space
(see \cite{Ru1} for similar spaces)
$$
BV^{(1)}=\{ \varphi : \real \to \complex\, ,
\supp(\varphi)\subset (-\infty, b]\, ,
\varphi' \in BV \}\, ,
$$
for the norm
$\|\varphi\|_{BV^{(1)}} =  \|\varphi'\|_{BV}$.
We have the following key lemma:

\begin{lemma}\label{spec}
The spectrum of $\LL_{0}$ on $BV^{(1)}$  and that of
$\LL_{1}$ on $BV$ coincide. 
In particular, 
the eigenvalues of modulus $> \lambda$ of the two
operators are
in bijection. 
\end{lemma}

\begin{proof}
By construction $\varphi \mapsto\varphi'$ is
a Banach space isomorphism between $BV^{(1)}$ and
$BV^{(0)}$. The Leibniz formula and the
chain rule imply that for any $\varphi \in BV^{(1)}$ 
\begin{equation}\label{keye}
(\LL_{0} \varphi)'=\LL_{1} (\varphi')\, .
\end{equation}
Indeed, the singular term in the Leibniz
formula (corresponding to the
derivative of  $\chi$, which is a dirac mass at
$c_1$) vanishes, because
$\psi_{+}(c_1)=\psi_{-}(c_1)=c$ and
$\varphi(c)-\varphi(c)=0$. 

That is,
the operators $\LL_{0}$ and $\LL_{1}$ are conjugated, and $\LL_{0}$
on $BV^{(1)}$
inherits the  spectral properties of $\LL_{1}$ on $BV$, as claimed.
\end{proof}

Lemma ~\ref{spec} implies that the spectral radius of $\LL_{0}$ on $BV^{(1)}$ is equal to $1$. 
The fixed vector is $R_0$, where
we define for $x\in \real$
\begin{equation}
R_0(x):=-1+\int_{-\infty}^x \rho_0(u) du\, .
\end{equation}
By construction, $R_0$ is Lipschitz, strictly increasing
on $[c_2, c_1]$, and constant outside of this
interval ($\equiv -1$ to the left and $\equiv 0$ to the right).
In addition, $R_0'$ coincides with $\rho_0$ on each continuity
point of $\rho_0$,  so that
$R_0'\sim \rho_0$.
The fixed vector of $\LL_0^*$ is $\nu(\varphi)=\varphi(b_0)-\varphi(a_0)$
with  $b_0=\psi_-(a_0)$.
Indeed, $\LL_0 \varphi(b)=0$ and $\LL_0 \varphi(a_0)=
\varphi(a_0)-\varphi(b_0)$. 
Since $b_0\ge b$ (otherwise we would have $\psi_- (a_0)=b_0> b$, a contradiction)
we have $\varphi(b_0)=\varphi(b)$.

\section{The susceptibility function and the decomposition $\rho_0=\rho_s+\rho_r$}

If $K$ is a compact interval we let $C^1(K)$ denote the set of functions on $K$
which extend to $C^1$ functions in an open neighbourhood of $K$.
The susceptibility function \cite{Ru3} associated to $f$ as above, $\varphi\in C^1([a_0,b])$, 
and the perturbation
$f_t=f+tX$, is  defined to be the formal power series
\begin{align}\label{susceptibility}
\Psi(z)&=\sum_{n=0}^\infty
\int z^n X (y) \rho_0(y) \frac {\partial} {\partial y} \varphi (f^n (y))\, dy \\
\nonumber &=\sum_{n=0}^\infty
\int z^n \LL^n_0(\rho_0 X)(x) \varphi' (x)\, dx\, .
\end{align}

The  expressions (\ref{susceptibility})
evaluated at $z=1$ may be obtained by formally
differentiating (\cite{Ru1ph}, see also  Appendix~\ref{susc} below)
the map 
\begin{equation}\label{srbmap}
\RR: t \mapsto \int \varphi(x) \rho_t(x)\, dx
\end{equation}
at $t=0$, when $\varphi$ is at least $C^1$.

\begin{proposition}\label{naive}
The power series $\Psi(z)$ extends to a holomorphic function in the open unit
disc, and in this disc we have 
$$
\Psi(z)=
\int (\id-z \LL_0)^{-1} (X  \rho_0) (y) \, \varphi' (x)\, dx \, .
$$
\end{proposition}

\begin{proof} The spectral properties
of $\LL_0$ on $BV$ (Lemma~\ref{spBV}) imply that for each $\delta > 0$ there is
$C$ so that $\|\LL_0^n \|_{BV} \le C(1+\delta)^n$,
so that $\Psi(z)$ is holomorphic in the open unit
disc. 
\end{proof}

\begin{remark}
Ruelle \cite{Ru3} studied $\Psi(z)$ for real-analytic multimodal maps
$f$ conjugated to a Chebyshev polynomial (e.g. the ``full" quadratic map
$2-x^2$ on $[-2,2]$).
In this nonuniformly expanding analytic setting, the susceptibility function is not holomorphic
in the unit disc: It is meromorphic in the complex
plane but has poles of modulus $< 1$.
(See also \cite{JR} for generalisations to
other real-analytic maps with preperiodic critical points,
and see \cite{BJR} for determinants giving the locations
of the poles when the dynamics is polynomial.) 
The study of real analytic non uniformly hyperbolic 
interval map with non preperiodic, but nonrecurrent, critical point 
is in progress \cite{Ruwip}.
\end{remark}

In order to analyse further $\Psi(z)$, let us next decompose
the invariant density $\rho_0$ into a singular and
a regular part:
Any function $\varphi: \real \to \complex$ of bounded variation, with regular discontinuities, can be uniquely
decomposed as $\varphi=\varphi_s+\varphi_r$, where 
the regular term $\varphi_r$ is continuous and of bounded variation
(with $\var (\varphi_r)\le \| \varphi\|_{BV}$),
while the singular (or ``saltus") term $\varphi_s$ is a
sum of jumps 
$$
\varphi_s = \sum_{u \in \SS} s_u H_u\, ,
$$
where $\SS$ is an at most countable set, $H_u(x)=-1$ if $x < u$,
$H_u(x)=0$ if $x > u$ and $H_u(u)=-1/2$, and
the $s_u$ are nonzero complex numbers so that
$\var (\varphi_s)=\sum_u |s_u| \le \| \varphi\|_{BV}$.
(See \cite{RN}, noting that our assumption that the discontinuities
of $\varphi$ are regular gives the above formulation.)
In the case when $\varphi$ is the invariant
density of a piecewise smooth and expanding interval
map, we have the following additional smoothness
of the regular term (this observation, which was inspired by
the analogous statement for nonuniformly expanding maps \cite{Ruwip},
seems new):

\begin{proposition}\label{decc}
Consider the decomposition $\rho_0=\rho_s+\rho_r$ of the invariant
density $\rho_0 \in BV$. Then $\rho_r \in BV^{(1)}$.
\end{proposition}

\begin{proof}
We shall use the following easy remark: If 
$a_0=x_0< x_1<\cdots <x_m=b$ for $m\ge 2$ and
$\varphi(x)=\varphi_i(x)$ for $x_{i-1}\le x \le x_i$,
with $\varphi_i$
extending to a $C^1$  function in a neighbourhood
of $[x_{i-1}, x_i]$, for $i=1, \ldots, m$,
then if $\varphi$ is supported in $[a_0,b]$ we have
$\varphi \in BV$ with 
\begin{equation}\label{bvbd}
\| \varphi \|_{BV}
\le (b-a_0) \sup_{i=1,\ldots, m} \sup_{[x_{i-1},x_i]}
|\varphi_i'| + \sum_{i=1}^{m-1}|\varphi_i(x_i)- \varphi_{i+1}(x_i)|
+ |\varphi_1(x_0)| + |\varphi_2(x_m)|\, .
\end{equation}

In this proof we write $\rho$ instead of $\rho_0$.
We know that if $\varphi_0\in BV$ is such that
$\int_{a_0}^b \varphi_0\, dx=1$ 
then $\rho=\lim_{n \to \infty} \rho^{(n)}$
with $\rho^{(n)}=\LL^n_1 (\varphi_0)$, the limit
being in the $BV$ topology.
We can assume in addition that $\varphi_0$ is $C^2$ and nonnegative. 
Decomposing $\rho^{(n)}=\rho^{(n)}_{s}+\rho^{(n)}_{r}$, we have
on the one hand that $\rho^{(n)}_{s}$ is a sum of jumps along
$c_j$ for $1\le j \le n$.
On the other hand, by the  remark in the beginning of the
proof, $\rho^{(n)}_{ r}$ is an element of $BV^{(1)}$.
We may estimate the $BV$ norm of $\Delta_n=(\rho^{(n)}_{ r})'$
as follows: First note that $\Delta_{n}$ extends to a $C^1$ function 
in a neighbourhood of $x$ 
if $x \notin\{ c_j\, ,1\le j \le n\}$. 
Next, we shall show by an easy distortion estimate  that there is $C$ (depending on $f$
and on the $C^2$ norm of $\varphi_0$) so that
\begin{equation}\label{est1}
|\Delta_n'(x)| \le C 
\, , \quad \forall n\, , \forall  x \notin\{ c_j \, ,1\le j \le n\}\, .
\end{equation}
Indeed, note  that if $x=f^n(y)$ with $x \notin\{ c_j \, ,1\le j \le n\}$
(so that $f^k(y) \ne c$
for $0\le k \le n-1$)
\begin{equation}\label{first}
\frac {d}
{dx}\frac{ 1}{ (f^n(y))'}=-
\sum_{k=0} ^{n-1} 
\frac{f''(f^k(y))}{f'(f^k(y))}
\frac{1}{(f^{n-k-1})'(f^{k+1}(y))}
\frac{1}{(f^n(y))'} \, .
\end{equation}
Since $\sup_{w\ne c} |f''(w)|/|f'(w)|\le C_0$ and $\sum_{k=0}^{n-1} |(f^{n-k-1})'(y_k)|^{-1}$
is bounded by a geometric series, uniformly in
$\{y_k \mid f^\ell(y_k) \ne c\, ,
\, 0\le \ell \le n-k-2\}$, we get 
$$|\Delta_n(x)| \le \widetilde C \LL^n_1 (\varphi_0) (x)+ \lambda^n \LL_1^n(|\varphi_0'|(x)\, ,
$$
where $\lambda \in (\sup_{x\ne c} |f'(x)|^{-1}, 1)$.
(We have not detailed the contribution of the terms where $\varphi_0$ has
been differentiated.)
The claim (\ref{est1}) follows from differentiating
the right-hand-side of (\ref{first}) with respect to $x$, and using  that
$\sup_{w \ne c} |f'''(w)|/|f'(w)| \le C_1$ and
$\sup_n \sup_x | \LL_1^n(\phi)|< \infty$ for all bounded $\phi$.

To conclude our analysis of the $BV$ norm of $\Delta_n$, we must consider
$x\in\{ c_j \, ,1\le j \le n\}$
and estimate
$|\lim_{w \uparrow x}\Delta_n(w)- \lim_{z\downarrow x}\Delta_n(z)|$.
The jump between the left and right limits corresponds
to the discrepancy between the sets $f^{-n}(w)$ and $f^{-n}(z)$,
i.e., it is of the same type as 
$|\lim_{w \uparrow x}\rho^{(n)}(w)- \lim_{z\downarrow x}\rho^{(n)} (z)|$,
with the difference that $1/|(f^{n})'(y)|$ or $\varphi_0(y)$ 
(for $f^n(y)=x$) are replaced by
their derivatives with respect to $x$. 
We find for all $n$ and all  $x \in\{ c_j \, ,1\le j \le n\}$
\begin{equation}\label{est2}
|\lim_{w \uparrow x}(\Delta_n)(w)- \lim_{z\downarrow x}(\Delta_n) (z)| \le \tilde C |\lim_{w \uparrow x}\rho^{(n)}(w)- \lim_{z\downarrow x}\rho^{(n)} (z)|\, .
\end{equation}
Thus, there is $\tilde C$ so that for all $n$
\begin{equation}\label{est3}
\sum_{  x \in\{ c_j \, ,1\le j \le n\}} |\lim_{w \uparrow x}(\Delta_n)(w)- \lim_{z\downarrow x}(\Delta_n) (z)| 
\le \tilde C \var (\LL^n_1(\varphi_0)) + \lambda^n \var (\LL_1^n(|\varphi_0'|)) \, .
\end{equation}

By the Lasota-Yorke estimates (see e.g. (\ref{LYu})) on
$\LL_1^n$, (\ref{est1}) and (\ref{est3}), together with
(\ref{bvbd}) imply that there is $\widehat C$ so that
$\| \Delta_n\|_{BV} \le \widehat C$ for all $n$.
Applying Helly's selection theorem, a subsequence
$\Delta_{n_k}$ converges pointwise and in $\LL^1(\mbox{Lebesgue})$
to some $\Delta \in BV$. Similar arguments show that $\rho^{(n_k)}_{r}$
and $\rho^{(n_k)}_{s}$ converge to $\hat \rho_r$ and $\hat \rho_s$, respectively (maybe
restricting further the subsequence). It follows that 
$\int \hat \rho_r \psi' \, dx=-\int \Delta \psi\, dx$
for all $C^1$ functions $\psi$, i.e. $\Delta=\hat \rho_r'$.
By construction we have
$\rho=\hat \rho_s+\hat \rho_r$, with
$\hat \rho_r\in BV^{(1)} \subset BV \cap C^0$,
and $\hat \rho_s$ a sum of jumps along the (at most countable)
postscritical orbit. By uniqueness of the decomposition
$\rho=\rho_s+\rho_r$, we have proved the lemma.
\end{proof}

We  may now consider the contribution  to $\Psi(z)$ of the
regular term in the decomposition from Proposition ~\ref{decc}:

\begin{lemma}\label{regg}
If $\varphi \in C^1([a_0,b])$ then
$$
\int (\id-z \LL_0)^{-1} (X  \rho_r) (x) \, \varphi' (x)\, dx \, 
$$
extends to a meromorphic function in a disc of radius strictly larger
than $1$, with only singularity in the closed unit disc
an at most simple pole at $z=1$.
The residue of this pole is
$
X(a_0)\rho_r(a_0) \bigl ( \int_{a_0}^b \varphi \rho_0\, dx -\varphi(a_0)
\bigr )
$.
\end{lemma}

\begin{proof}
The spectral properties of $\LL_0$ on $BV^{(1)}$ (Lemma~\ref{spec})
imply that $(\id-z\LL_0)^{-1}(X\rho_r)$ depends
meromorphically on $z$ in a disc of radius strictly larger than $1$,
where its only possible singularity in the closed unit disc
is a simple pole at $z=1$, with residue
$(X(b_0)\rho_r(b_0)-X(a_0)\rho_r(a_0)) R_0(x)$.
Since $\rho_r$ is continuous and supported in $(-\infty,c_1]\subset
(-\infty,b_0]$ we have $\rho_r(b_0)=0$. 
To finish, integrate $\int_{a_0}^b \varphi' R_0 \, dx$
by parts and use $R_0(b)=0$ and $R_0(a_0)=-1$.
 \end{proof}

Clearly, $(\id-z\LL_0)^{-1}(X\rho_s)=
\sum_{n=0}^\infty z^n\LL_0^n(X\rho_s)$ is an element of $BV$ which depends 
holomorphically on $z$ in the open unit disc. We will
be able to say much more about this expression if $c$ is preperiodic,
in Section~\ref{markov}.
If $c$ is not preperiodic, the situation is not as transparent, but some
results are collected in Section~\ref{nonM}.
In view of Sections~\ref{nonM}--\ref{markov}, we introduce further notation.

If $c$ is preperiodic, i.e. $f^{n_0}(c)$ has minimal
period $n_1\ge 1$ (with $n_0 \ge 2$ minimal), we set $N=n_0+n_1-1\ge 2$, otherwise we
put $N=\infty$. By definition of
the saltus, we  have
\begin{equation}\label{defsalt}
\rho_s(x)=\sum_{n=1}^N s_n H_{c_n}(x)\, ,
\end{equation}
with
$
s_n=\lim_{y \downarrow c_n} \rho(y)-\lim_{x \uparrow c_n} \rho (x)
$.

We next define the weighted total jump of $f$:
\begin{align}\label{saltusid}
\JJ(f,X)&=\sum_{n=1}^N s_n X(c_n) \, .
\end{align}
We put $\JJ(f)=\JJ(f,1)$.
Note that  
$$\JJ(f)=-\rho_r(b_0)+\rho_r(a_0)=\rho_r(a_0)\, .
$$

\begin{remark}
If $f$ is a tent-map, i.e. $|f'(x)|$ (for $x\ne c$)
is constant, then it is easy to see that $\rho=\rho_s$ is purely
a saltus function (for example use  $\rho=\lim_{n\to \infty}
\LL_1^n (\varphi_0)$, with $\varphi_0$ the normalised
characteristic function of $[c_2, c_1]$, the limit being in the variation norm). In particular,
we get that $\JJ(f)=\JJ(f,1)=0$ for all tent-maps.
\end{remark}


\section{The susceptibility function in the non-Markov case}
\label{nonM}

In this section we assume (i)--(iv) and  that $c$ is not preperiodic
(i.e. for every $n\ge 1$, the point $f^n(c)$ is not periodic;
in other words, there does not exist a finite Markov partition
for $f$). We can suppose
without further restricting generality that $f(c)< b$ and $\min(f(a), f(b))> a> a_0$.
We start with a preparatory lemma:

\begin{lemma}\label{candidate}
Assume that $c$ is not preperiodic.

If  $\JJ(f)=0$  then the function
$
\tilde \rho_s=\sum_{j=1} ^\infty H_{c_j} \sum_{k=1}^j s_k 
$
is of bounded variation and
satisfies $(\id-\LL_0)\tilde \rho_s=\rho_s$.

If  $\JJ(f,X)=0$ then,  setting  $\delta_{c_j}$ to be the dirac mass at
$c_j$, the measure 
$
\mu_s=\sum_{j=1} ^\infty \delta_{c_j}  \sum_{k=1}^j s_k X(c_k)
$
is bounded and satisfies
$
(\id- f_*) \mu_s  =X \rho_s'
$.
\end{lemma}

\begin{remark}
We do {\it not} claim that when $\JJ(f)=0$
the sum $\sum_{n=0}^\infty \LL^n_0 \rho_s$ converges to $(\id-\LL_0)^{-1}\rho_s=\tilde \rho_s$
or that $(\id -z\LL_0)^{-1}\rho_s$ converges to
$(\id-\LL_0)^{-1}\rho_s$ as $z \to 1$ (even within $[0,1]$),
and we do not claim the parallel statements 
about $(\id - f_*)^{-1}(X \rho'_s)=\mu_s$ when $\JJ(f,X)=0$.
\end{remark}

\begin{remark}
For any complex number $\kappa$ we have
$(\id-\LL_0)(\tilde \rho_s+\kappa R_0)=\rho_s$
and $(\id- f_*) (\mu_s+\kappa \rho_0)  =X \rho_s'$.
\end{remark}

Our  result in this case is:

\begin{proposition}\label{prop1}
Assume that $c$ is not preperiodic
and let $\varphi \in C^1([a_0,b])$.

For $|z|< 1$ we have
\begin{align}\label{chsusceptibility}
\Psi(z)
&=-\sum_{j=1}^\infty \varphi(c_j)  \sum_{k=1}^j z^{j-k} s_k X(c_k) 
-
\int (\id- z\LL_1)^{-1} (X' \rho_s+(X\rho_r)')(x) \varphi (x)\, dx  \, .
\end{align}
The second term above extends to a meromorphic function in a disc
of radius strictly larger than $1$, with only singularity an at most simple
pole at $z=1$, with residue
$
\JJ(f, X) \int_{a_0}^b \varphi \rho_0 \, dx 
$.

If $\JJ(f,X)=0$ then the following is a well-defined complex number:
\begin{align}\label{psi1}
\Psi_1
&=-\sum_{j=1}^\infty \varphi(c_j)  \sum_{k=1}^j  s_k X(c_k)
-
\int (\id- \LL_1)^{-1} (X' \rho_s+(X\rho_r)')(x) \varphi (x)\, dx \, .
\end{align}
\end{proposition}

\begin{remark}\label{smania}
There 
exists a unique function  
$\alpha$ on the postscritical orbit
so that $X(c_{k+1})=\alpha(c_{k+1})-f'(c_k)\alpha(c_k)$
for $k\ge 1$: 
set $$\alpha(c_k)=-\sum_{j\ge 0} X(c_{k+1+j})/(f^{j+1})'(c_k)\, .
$$
(See e.g. \cite[Proof of Thm 1]{Sm} for the relevance of this
``twisted cohomology equation", in view of Conjecture A:
The possibility to extend $\alpha$ ``smoothly" to  $I$ is related
to the existence of a topological conjugacy between $f$ and $f_t$.)
Since $s_k=f'(c_k)s_{k+1}$ for all
$k \ge 1$, and since $s_1\ne 0$, our condition $\JJ(f,X)=0$ is equivalent
to requiring that $X(c_1)-\alpha(c_1)=0$.
\end{remark}

In view of Lemma~\ref{candidate}, slightly abusing notation, we may write when 
$c$ is not preperiodic,  and $\JJ(f,X)=0$ 
$$
 \Psi_1
=
-\int (\id- \LL_1)^{-1} ((X \rho_0)')(x) \varphi(x) \, dx\, .
$$
If, in addition, $X \equiv 1$, we may 
also write
\begin{align}\label{chsusceptibility2}
&\Psi_1=
\int (\id - \LL_0)^{-1} (\rho_0)(x) \varphi' (x)\, dx  \, .
\end{align}

The orbit of $c$ is expected to be
``generically" dense, so that  both conditions
``$\varphi(c_k)=0$ for all $k$ and $\int \varphi \rho_0\, dx=0$"
and   ``$X(c_k)=0$ for all $k$" are very strong.
\footnote{They are satisfied for nontrivial  $X$ e.g.
if $c$ is not recurrent.}
However, we point out that 
either condition implies that $\Psi(z)$ extends
holomorphically to a disc of radius larger than $1$, with $\Psi(1)=\Psi_1$.

The relationship between $\Psi(z)$ and
$\Psi_1$ (when $\JJ(f,X)=0$) is unclear
for general $\varphi$ and $X$.
(See Remark~\ref{corr}. See however Appendix \ref{regs}
for an alternative -- perhaps artificial -- 
susceptibility function, which can be related to $\Psi_1$.)
If $\JJ(f,X)\ne 0$, it seems unlikely that a replacement 
for $\Psi_1$ would exist. (See also Appendix \ref{regs}.)

\medskip

We now prove Lemma~\ref{candidate}.

\begin{proof}
Note first that if $c$ is not preperiodic then, since $\rho=\lim_{j \to \infty} \LL_1^j (\varphi_0)$ 
(for $\varphi_0$ as in the proof of Proposition~\ref{decc}) and
the convergence is exponentially rapid in the $BV$ norm,
there are $\xi < 1$ and  $C\ge 1$ so that 
\begin{equation}\label{notper}
\sum_{k\ge j+1}  |s_k| \le \var (\rho^{(j)} - \rho) \le C \xi^j \, , \forall j \, .
\end{equation}

Then apply (\ref{notper}) and the assumption  $\JJ(f)=\sum_{k=1}^\infty s_k=0$,
to get
\begin{equation}\label{zerojump}
|\sum_{k=1}^j s_k|=|-\sum_{k\ge j+1} s_k | \le C \xi^j\, , \forall j\, .
\end{equation}
Observe next  that
$\LL_0 (H_{c_j})= H_{c_{j+1}}$ for all $j \ge 1$.
Finally, use  $\sup |H_{c_j}|\le 1$ for all $j$ and (\ref{defsalt}). 

For the second claim, use also  $\sup |X|\le 1$, $f_*(\delta_{c_j})=\delta_{c_{j+1}}$, 
that $\JJ(f,X)=0$,
implies
\begin{equation}\label{zerojumpX}
|\sum_{k=1}^j s_k X(c_k) | \le C \xi^j\, , \forall j\, .
\end{equation}
and that ($\rho_s'$ is a distribution of
order $0$ and $X$ is continuous)
\begin{equation}\label{Xr'}
X \rho_s'=\sum_{j=1}^\infty X(c_j) s_j  \delta_{c_j}\, .
\end{equation}
\end{proof}

We next show Proposition ~\ref{prop1}:

\begin{proof}
Write $\rho$ for $\rho_0$ and consider the
decomposition $\rho=\rho_s+\rho_r$. 
We have $\rho(b)=\rho(b_0)=\rho(a_0)=0$,  $\rho_r(b)=0$,
and $\rho_s$ is continuous at $a_0$, $b_0$, and $b$ with $\rho_s(b_0)=\rho_s(b)=0$.
We may integrate
by parts,  and get from the Leibniz formula
(recall Lemma ~\ref{spec}, and note that $(\LL_0(\psi))'=f_*\psi'$,
for $\psi\in BV$) for $|z|< 1$ that
\begin{equation}\label{Leib}
\int_{a_0}^b  \LL^n_0(\rho X)(x) \varphi' (x)\, dx=
-  \int \varphi f^n_* (X \rho_s')-
\int \LL^n_1(X'\rho_s+(X\rho_r )')(x) \varphi (x)\, dx\, .
\end{equation}
(There are no boundary terms in the Stieltjes integration
by parts because $\LL_0^n(\rho X)$ is continuous and vanishes at $b$, 
$\rho(b_0)=0$ and $\rho(a_0)=0$.)
It follows that for $|z|< 1$
\begin{align}
\Psi(z)&=-\sum_{n=0}^\infty z^n \biggl ( \int \varphi f^n_* (X \rho_s')+
\int \LL^n_1(X'\rho_s+(X\rho_r )')(x) \varphi (x)\, dx\biggr )\, .
\end{align}
The proof of Lemma~\ref{regg} applies to $\LL_1$ on  $BV$
and allows us to control the terms associated to
$(X\rho_r)'$ and $X' \rho_s$. Since $\JJ(f,X)=\int_\real X \rho'_s$, 
the residue of the possible pole
at $z=1$ is, using Stieltjes integration by parts,
\begin{align}\label{stj2}
\nonumber -(\int_\real X' \rho_s\, dx + \int_\real (X \rho_r)'\, dx) &\int_{a_0}^b \varphi \rho_0 \, dx
=-(\int_\real X' \rho\, dx + \int_\real X \rho_r'\, dx) \int_{a_0}^b \varphi \rho_0 \, dx\\
\nonumber &=(\int_\real X \rho' - \int_\real X \rho_r'\, dx)\int_{a_0}^b \varphi \rho_0 \, dx\\
&=\JJ(f, X) \int_{a_0}^b \varphi \rho_0 \, dx \, .
\end{align}
\footnote{In the case $X=1$, recall that
$\JJ(f,1)= \rho_r(a_0)$, and note that
$\int\LL_0^n(\rho_s)\varphi' \, dx=-\int \varphi f^n_* (\rho_s')
+\varphi(a_0)\rho_r(a_0)$.}
On the other hand, 
we get by (\ref{notper}), (\ref{Xr'}), and  since $\sup |H_{c_j}|\le 1$ for all $j$, 
that
for each $|z|< 1$
\begin{align}\label{exchange}
\sum_{n=0}^\infty z^n f^n_* (X \rho_s')  
&=\sum_{n=0}^\infty z^n \sum_{k=1}^\infty  s_{k} X(c_{k}) \delta_{c_{k+n}}\\
\nonumber &= \sum_{j=1}^\infty   \delta_{c_j} \sum_{k=1}^j z^{j-k} X(c_k) s_k \, .
\end{align}
We have proved (\ref{chsusceptibility}) in the open disc of radius
$1$.
The fact that $\Psi_1$ is well-defined follows from Lemma~\ref{candidate}
and our assumption that $\JJ(f,X)=0$ which implies (\ref{zerojumpX}).

\end{proof}

\begin{remark}\label{corr}
In spite of Lemma~\ref{candidate}, we are not in a position to apply
Fubini's theorem in (\ref{exchange}) at $z=1$. It seems unlikely that the sum $\sum_{n=0}^\infty  \sum_{k=1}^\infty s_k  X(c_{k})\delta_{c_{k+n}}$
converges in the usual sense to $\mu_s$, and it is unclear whether
$\mu_s$ could be interpreted as a classical (e.g. Norlund or Abelian) limit of this sum. 
\end{remark}


\section{The susceptibility function in the Markov case}
\label{markov}

Assume in this section, in addition to (i)--(iv), that $c$ is preperiodic,
i.e. there exist $n_0\ge 2$ and $n_1\ge 1$
so that $c_{n_0}$ is periodic of minimal period $n_1$
(we take $n_0$ minimal for this property).
In this Markov case, we have the following result:

\begin{theorem}\label{M1}
Assume that $c$ is preperiodic. 
Let $\varphi \in C^1([a_0,b])$.
Then $\Psi(z)$ admits a meromorphic extension to a disc of radius 
$>1$.  The poles of $\Psi(z)$
in the closed unit disc are at most simple poles at the
$n_1$th roots of unity.

Assume  either 
(a) $\varphi(c_k)=0$ for all
$k \ge n_0$ and $\int \varphi \rho_0 \, dx=0$,
or (b) $X(c_k)=0$ for all $k \ge 1$, then  the residues of all 
poles of modulus one  of $\Psi(z)$ vanish, and  
$$
\Psi(1)
=\lim_{k \to \infty}\sum_{n=0}^k
\int  \LL^n_0(\rho_0 X)(x) \varphi' (x)\, dx \, .
$$
\end{theorem}

We next exhibit other sufficient conditions for the residues of the poles
of $\Psi(z)$ on the unit circle
to vanish.  For this, we introduce 
$\JJ^{1,n_0}_{n_0}=\JJ^{1,n_0}_{n_0}(f,X)=\JJ(f,X)=\sum_{k=1}^{n_0} X(c_k) s_k$, and, if
$n_1 \ge 2$, the following sums of jumps for $m=n_0, \ldots,n_0+ n_1-1$:
\begin{align*}
\JJ^{n_1,n_0}_{m}&=\JJ^{n_1,n_0}_{m}(f,X)=\sum_{
\substack{1\le k\le n_0+n_1-1:\\
\exists \ell \ge 0: k+n_0-1-\ell n_1=m
}} X(c_k) s_{k} \, .
\end{align*}

\begin{theorem}\label{M2}
Assume that $c$ is preperiodic. 
Let $\varphi\in C^1([a_0,b])$.

If $\JJ^{n_1,n_0}_m=0$ 
for $m=n_0, \ldots n_0+n_1-1$,
then  $\Psi(z)$ is holomorphic 
in a disc of radius strictly larger than one
with
$$
\Psi(1) =\lim_{z\to 1} \Psi(z)
\mbox{ and } \Psi(1)
=\lim_{k \to \infty}\sum_{n=0}^k
\int  \LL^n_0(X \rho_0)(x) \varphi' (x)\, dx\, .
$$

The residue of $\Psi(z)$ at $z=1$ is
$\JJ(f,X) (\int \varphi \rho_0 \, dx-
\frac{1}{n_1}\sum_{j=n_0}^{n_0+n_1-1} \varphi(c_j))$,
in particular, if $\JJ(f,X)=0$ then  $\Psi(z)$ is holomorphic at
$z=1$ with
$$
\Psi(1) =\lim_{z\in [0,1), z\to 1} \Psi(z)\, .
$$
\end{theorem}

We first prove Theorem~\ref{M1}:

\begin{proof}
Since
$\LL_0^n (\rho X)$ is continuous and vanishes 
at  $b$,   the
term associated to the  rightmost boundary in  the Stieltjes integration
by parts (\ref{Leib})  in the proof of Proposition ~\ref{prop1} vanishes.
If $c_2\ne a_0$ then $\rho$  vanishes and
is continuous at $a_0$ and $b_0$, so that the leftmost boundary term
from (\ref{Leib}) vanishes.
If $c_2=a_0$, this leftmost boundary term is in fact included
in the Stieltjes integral $-\int \varphi f^n_*(X\rho_s')$.

We consider $X\equiv 1$, the general case follows by
integration by parts as in (\ref{Leib}--\ref{stj2})
in the proof of Proposition ~\ref{prop1}
(recall in particular the residue $\JJ(f,X)\int \varphi \rho_0\, dx$), using
the remarks in the previous paragraph. 
By Lemma ~\ref{spec}, Proposition~\ref{decc}, Lemma~\ref{regg}
and (\ref{defsalt}) it suffices to consider $\LL_0$ acting on the finite-dimensional space
generated by $H_{c_k}$, for $1\le k \le n_0+n_1-1$.
We have 
$$
\begin{cases}
\LL_0 H_{c_j}=H_{c_{j+1}}\, , & j < n_0+n_1-1\,,\\
\LL_0 H_{c_{n_0+n_1-1}}=H_{c_{n_0}}\, .
\end{cases}
$$
The $(n_0+n_1-1) \times (n_0+n_1-1)$ matrix $L$ associated to the above linear operator
is such that $L^{n_1}$ is in lower triangular form, with zeroes in
the first $n_0-1$ diagonal elements and
with an $n_1\times n_1$ identity block in the $n_1$ last rows and columns.
It follows that  $(\id-z \LL_0) ^{-1}(\rho_0)$ 
extends meromorphically to a disc of radius strictly larger than
$1$, whose singularities 
on the unit circle are at most simple poles at the
$n_1$th roots of unity. 
To show the claim on the vanishing of the
residues, we integrate $\int_{a_0}^b \LL_0^n(\rho_s)\varphi' \, dx$
by parts: it suffices to consider
the boundary terms since our assumption
$\varphi(c_j)=0$ for all $j\ge n_0$ guarantees that the poles
corresponding to the eigenvalues of $L$ have zero residue.
If $a_0\ne c_2$ then the boundary term gives a residue $-\varphi(a_0)\rho_s(a_0)$
for the pole at $z=1$,
which, summed with
the residue from Lemma~\ref{regg} 
gives $\JJ(f,1) \int \varphi \rho_0 \, dx$ (using $\rho_s(a_0)+\rho_r(a_0)=0$
and $\JJ(f,1)=\rho_r(a_0)$).
If $a_0=c_2$, the boundary term gives rise to the multiple of $\varphi(c_2)$ 
which appears in the contribution of the spectrum of $L$, and
Lemma~\ref{regg} gives
$\JJ(f,1)(\int \varphi \rho_0\, dx-\varphi(a_0))$.
\end{proof}

We now prove Theorem ~\ref{M2}:

\begin{proof}
Again, we consider $X\equiv 1$, and the general case follows by
integration by parts.
If $\JJ^{n_1,n_0}_m=0$ 
for $m=n_0, \ldots n_0+n_1-1$,
then $\LL_0^{n_0-1}(\rho_s)$ vanishes. It follows
(recall Lemma~\ref{regg}, the residue there
vanishes if $X=1$ since $\rho_r(a_0)=\JJ(f,1)$) that $\Psi(z)$ is holomorphic in a disc 
of radius strictly larger than one.

If $\JJ(f)=0$ then we claim that the spectral projector $\Pi$ associated to the eigenvalue
$1$ of the matrix $L$ introduced in the proof of Theorem ~\ref{M1}
satisfies $\Pi(\rho_s)=0$ (this gives the second statement of the theorem).
To show the claim note that the fixed vector for $L$
is $v=(v_j)$ with $v_j=0$ for $j\le n_0-1$
and $v_j=1$ for $n_0\le  j \le n_0+n_1-1$,
and that $u=(1, \ldots, 1)$ is a left fixed vector
for $L$. The projector $\Pi$ is just $\Pi(w)= \frac{\langle u,w\rangle }
{\langle u,v\rangle }v$, and $\Pi(\rho_s)=0$
follows from $\JJ(f,1)=0$.
\end{proof}


\section{Non differentiability of the SRB measure}
\label{countt}

In this section, we present examples
\footnote{The example in Theorem~\ref{counter} and Remark \ref{AAA}
are due to A. Avila.
D. Dolgopyat told me several years ago that he believed
the SRB measure was not a Lipschitz function of the dynamics
in the present setting, and he may have been aware of similar examples.
After this paper was written, we learned about
\cite{Er} which, although mostly nonrigorous, indicated that 
$\RR(t)$ should not be expected to be Lipschitz, and
C. Liverani brought to our attention Mazzolena's  \cite{MM}
detailed analysis of families of maps for which $\RR(t)$
is not Lipschitz.} 
of perturbations
$f+tX\circ f$ of maps $f$ satisfying (i)--(iv), 
so that $f$ has a
preperiodic critical point, and at which $t\mapsto \int \varphi \rho_t \, dx$
fails to be Lipschitz at $t=0$
for a well-chosen smooth observable $\varphi$.
(In view of (\ref{claim'}), we shall see that
the examples are ``as bad as possible.")

Recall that we call tent-map a map $f$ satisfying (i)--(iv) and so that
$|f'(x)|$ is constant for $x \ne c$.
For $1< \lambda \le 2$ we let
$g_\lambda$ be the tent-map of slopes $\pm \lambda$
on $[0,1]$, i.e.,
$g_\lambda(x)=\lambda x$ for $x \in [0,1/2]$, and $g_\lambda(x)=\lambda-\lambda x$ for 
$x \in [1/2,1]$. We put $c_n(\lambda)=g_\lambda^n (1/2)$ for $n \ge 1$.

We first present the simplest possible counter-example:

\begin{theorem}\label{counter}
There exists  a $C^1$ function $\varphi$, with $\varphi(0)=\varphi(1)=0$,   a
sequence $\lambda_k\in (1,2)$ with  $\lim_{k \to \infty} \lambda_k =2$, 
so that  $c_{k+2}(\lambda_k)$ is a fixed point of $g_{\lambda_k}$, and 
a constant $C > 0$ so that
$$
\int \varphi \gamma_{\lambda_k} \, dx - \int \varphi \gamma_2 \, dx
\ge C k (2-\lambda_k)\, ,
\quad \forall k \, ,
$$
with $\int \varphi \gamma_2 \, dx=1$, where $\gamma_{\lambda_k}$ is the
invariant density of $g_{\lambda_k}$.
\end{theorem}

(In fact we have $\varphi(c_n(\lambda_k))=0$ for all $k\ge 1$ and $n\ge 1$
in Theorem~\ref{counter}.)

The theorem shows that the SRB measure cannot be
(one-sided) Lipschitz at $g_2$ for $\varphi$.
Since we can write $g_{\lambda_k}=g_2+t_k X\circ f$, with $t_k = \lambda_k-2$ and
$X$ as in \S2, 
with $X(0)=0$ (in fact, $X(x)=x$ for $x\in [0,1]$), and $\varphi(0)=\varphi(1)=0$, 
Theorem~\ref{M2}  applies to $f=g_2$, $X$,  and $\varphi$,
and, since $\int \varphi\rho_0\, dx\ne 0$ gives that $\Psi(z)$ is meromorphic in a disc of radius larger
than one with a simple pole at $z=1$ (the residue is
$\JJ(f,X)\int \varphi \rho_0\, dx$ with $\JJ(f,X)=X(c_1)s_1\ne 0$).
Note that Ruelle \cite{Ru3} proved that the susceptibility function
associated to the full quadratic map and any smooth $X$ and $\varphi$
has a vanishing residue at $z=1$. However, $\Psi(z)$ has a pole
strictly inside the unit disc in the setting of \cite{Ru3}.

Of course, the example in Theorem~\ref{counter} is a bit special since
$g_2$ is an ``extremal" tent-map. But it is not very difficult to
provide other examples of tent-maps with preperiodic critical points
at which the SRB measure is not a Lipschitz function of the dynamics.
Indeed, coding the postcritical orbit by the sequence $\Theta$,
with $\Theta_j=L$ if $c_j<1/2$ and $\Theta_j=R$ if $c_j > 1/2$, 
the code of $g_2$ is $RL^\infty$
(that is, $\Theta_1=R$, and $\Theta_j=L$ for
all $j\ge 2$), while the proof of Theorem~\ref{counter}
shows that the code of $g_{\lambda_k}$ is $\Theta_1=R$,
$\Theta_j=L$ for $2\le j \le k+1$ and $\Theta_j=R$ for $j \ge k+2$.
The following example corresponds to a similar perturbation,
starting from $\Theta=RLR^\infty$ (i.e., $g_{\sqrt 2}$),
and considering a sequence $g_{\nu_\ell}$, for $\ell \ge 6$ and
even, where $\nu_\ell$ is the unique parameter giving the code
\begin{equation}\label{defl}
\Theta_1=R\, , \, \Theta_2=L\, , \, \Theta_j=R \mbox{ for } 3\le j \le \ell-2\, ,
\, 
\Theta_{\ell-1}=L\, , \,  \Theta_j=R \mbox{ for } j \ge \ell \, .
\end{equation}
(In particular $c_\ell(\nu_\ell)$ is the fixed point of $g_{\nu_\ell}$.)

\begin{theorem}\label{counter2}
There exists  a $C^1$ function $\varphi$, with $\varphi(c_1(\sqrt 2))=\varphi(c_2(\sqrt 2))=
\varphi(c_3(\sqrt 2))=0$, and $\int \varphi \gamma_{\sqrt 2}\, dx =1$, a
sequence $\nu_\ell\in (\sqrt 2,2)$, with $\ell$ even
and  $\lim_{\ell \to \infty} \nu_\ell =\sqrt 2$, 
so that  $c_{\ell}(\nu_\ell)$ is a fixed point of $g_{\nu_\ell}$, and 
a constant $C > 0$ so that
$$
\int \varphi \gamma_{\sqrt 2}\, dx  -\int \varphi \gamma_{\nu_\ell} \, dx 
\ge C \ell (\nu_\ell-\sqrt 2)\, ,
\quad \forall \ell \, .
$$
\end{theorem}

(In fact we have $\varphi(c_n(\nu_\ell))=0$ for all even $\ell\ge 4$ and $n\ge 1$
in Theorem~\ref{counter2}.)

Theorem~\ref{M2} applies to the example in Theorem~\ref{counter2}
and gives that $\Psi(z)$ has a simple pole at $z=1$ with residue
$\JJ(f,X) \int \varphi \gamma_{\sqrt 2}\, dx \ne 0$.

\begin{remark}\label{AAA}
Although the combinatorics will be more complicated, a modification of
the proof of Theorems~\ref{counter} and ~\ref{counter2}
should be applicable \cite{AA} to all
preperiodic  tent-maps. This would give a dense
countable set of parameters $\Lambda_0$, and $C^1$ functions
$\varphi_{\lambda}$, where the  SRB measure $\lambda \mapsto
\RR(\lambda)=\int \varphi_{\lambda_0} \gamma_\lambda\, dx$ is
not Lipschitz at $\lambda_0$ if $\lambda_0\in \Lambda_0$, for which
$\Psi(z)$ is meromorphic at $z=1$ (by Theorem ~\ref{M2}).
If this construction is possible, a Baire argument \cite{AA}
would then imply that there is an uncountable set of parameters
$\Lambda_1$ where $\RR(\lambda)$ is not Lipschitz. This would
give rise to counterexamples which are non-Markov tent-maps
to which Proposition~\ref{prop1} applies
(with, presumably, $\JJ(f,X)\ne 0$).
\end{remark}

In view of the program sketched in the
previous remark, it would seem that the SRB measure of tent-maps is 
not often Lipschitz.

We next prove Theorem~\ref{counter}:

\begin{proof}
The fixed point of $g_\lambda$
is $x_\lambda=\lambda/(1+\lambda)>1/2$ and its preimage
in $[0,1/2]$ is $y_\lambda=1/(1+\lambda)$. Let
$z_\lambda=y_\lambda/\lambda$ be the preimage
of $y_\lambda$ in $[0,1/2]$. The critical value
is $c_1=\lambda/2 > 1/2$ (in this proof we write $c_j$ for $c_j(\lambda)$
whenever the meaning is obvious), which is mapped to $c_2=(2-\lambda)\lambda/2< 1/2$.
If $\lambda=2$ then $c_1=1$, $c_2=c_3=0$, and $\gamma_2$ is constant, equal to
$1$ on $[0,1]$.

If $1<\lambda < 2$, then $c_{k+1}=y_\lambda$ if $\lambda =\lambda_k= 2-C_k \lambda^{-k}$
with $C_k=2/(\lambda(1+\lambda))$, and $k \ge 1$.
The invariant density for such $\lambda_k$
is supported in $[c_2,c_1]$ and constant on
each 
$(c_{j+1}, c_{j+2})$ for $1\le j \le k$, 
with value $v_j > 0$,  and on $(c_{k+2},c_1)=(x_{\lambda_k},c_1)$, with value
$v_{k+1}> 0$. The fixed point equation for $\gamma_{\lambda_k}$
reads
$v_{k+1}=\lambda_k v_1$, $v_{j}+v_{k+1}=\lambda_k v_{j+1}$ for $1\le j \le k-1$,
and 
$2 v_{k}=\lambda_k v_{k+1}$. This implies that the sequence $j \mapsto v_j$ is strictly
increasing.
(Indeed, $v_{k+1}=2v_k/\lambda_k  > v_k$,
and proceed by decreasing induction, using
that $v_{j+1} > (v_j+v_{k+1})/2$
and $v_{k+1}> v_{j+1}$   to show that $v_{j+1} > v_j$ for $k-1\ge j \ge 1$.)
We take a nonnegative $C^1$ function $\varphi$ which is supported in 
$(2/3, 3/4)$, and thus in $(c_{k+2},c_1)$ for all large enough $k$. 
We  assume that $\int \varphi (x) \gamma_2\, dx=
\int \varphi(x)\, dx=1$. 
We next show that there is $D> 0$ so that
for all $k\ge 1$
\begin{equation}\label{toobig}
\int \varphi \gamma_{\lambda_k} \, dx 
\ge 1+ D k \lambda_k^{-k} \, ,
\end{equation}
and
this will end the proof of the theorem.

To show (\ref{toobig}), we use the fact that $\int \varphi \gamma_{\lambda_k} \, dx =
v_{k+1}\int_{c_{k+2}}^{c_1} \varphi(x)\, dx=v_{k+1}$.
To estimate $v_{k+1}$ we exploit $\int \gamma_{\lambda_k}\,  dx=1$:
This integral is equal to the difference 
$$
v_{k+1} (c_1-c_2)
-\sum_{j=1}^{k} (c_{j+2}-c_{j+1})(v_{k+1}-v_j)\, .
$$
We have $c_{j+2}-c_{j+1}= \lambda_k^{j-1} (c_3-c_2)$ with
$(c_3-c_2)\ge A \lambda_k^{-k}$ with $A$ independent of $k$, and
$(v_{k+1}-v_j)\ge (v_{j+1}-v_j)=v_{k+1} \lambda_k^ {k-j}(2-\lambda)/2\ge C_k\lambda_k^{-j}/2$
for $1\le j \le k-1$
\footnote{To get the equality, first use $v_m+ v_{k+1}=\lambda_k v_{m+1}$
at $m=j$ and $m=j+1$, repeat this  $k-j-1$ more times, and end by using
that $2 v_k = \lambda_k v_{k+1}$.},
so that $\sum_{j=1}^{k} (c_{j+2}-c_{j+1})(v_{k+1}-v_j)\ge E k \lambda_k^{-k-1}$
and
$$
1\le v_{k+1} (c_1-c_2)- Ek \lambda_k^{-k-1} \, ,
$$
which implies  $v_{k+1} \ge (1+ E k \lambda_k^{-k-1} )/(c_1-c_2)$.
Since $c_1-c_2 \le 1$  we  proved (\ref{toobig}).
\end{proof}

Finally, we show Theorem~\ref{counter2}:

\begin{proof}
Note that $\gamma_{\sqrt 2}$ is constant equal to
$u$ on $(c_2(\sqrt 2), c_3(\sqrt 2)$,
and  constant equal to $\sqrt 2 u$ on $(c_3(\sqrt 2),c_1(\sqrt 2))$, with $c_3(\sqrt 2)> 1/2$ the fixed point.
Putting
\begin{equation}\label{dd}
d=\bigl ( c_3(\sqrt 2)-c_2(\sqrt 2) + \sqrt 2(c_1(\sqrt 2)-c_3(\sqrt 2) \bigr ) \, ,
\end{equation}
the  normalisation condition is
\begin{equation}\label{noru}
du=1 \, .
\end{equation}

For $\ell \ge 6$ even, we define $\nu_\ell<\sqrt 2$ by (\ref{defl}).
Then $c_\ell(\nu_\ell)> c_3(\sqrt 2)$ is a fixed point
and the critical orbit of $g_{\nu_\ell}$ is ordered as follows
(in the remainder of this proof we write $c_m$ for $c_m(\nu_\ell)$
when the meaning is clear)
$$
c_2 < c_{\ell-1} < c_0 
< c_{\ell-3} < \cdots < c_5 <c_3 
<c_\ell 
< c_4 < c_6 < \cdots< c_{\ell-2} < c_1 \, .
$$
The invariant density of $g(\nu_\ell)$ is constant equal to $u_1=u_1(\nu_\ell)$ on $(c_2, c_{\ell-1})$,
constant equal to $u_2$ on $(c_{\ell-1}, c_{\ell-3})$,
constant equal to $u_j$ on $(c_{\ell-(2j-3)}, c_{\ell-(2j -1)})$
for $3\le j \le \ell/2-1$,
constant equal to $u_{\ell/2}$ on $(c_3, c_\ell)$,
constant equal to $u_{\ell/2+1}$ on $(c_\ell, c_4)$,
constant equal to $u_j$ on $(c_{2j-\ell}, c_{2j+2-\ell})$
for $\ell/2+2\le j \le \ell-2$,
constant equal to $u_{\ell-1}$ on $(c_{\ell-2}, c_1)$.
As $\ell \to \infty$ we have  that $c_1$ tends to
$c_1(\sqrt 2)$, that $c_2$ and $c_{\ell-1}$
tend to $c_2(\sqrt 2)$,
that $c_\ell$ 
tends to $c_3(\sqrt 2)$.
In particular, (\ref{claim'}) implies that $u_2 \to u$.

The fixed point equation for $\gamma_{\nu_\ell}$
implies that $u_{\ell-1}=\nu_\ell u_1$,
$u_{\ell-2}=\nu_\ell u_2$ and $2 u_2 = \nu_\ell u_{\ell-1}$
(thus, $u_{\ell-1}=\nu_\ell u_2/2$ tends to $\sqrt 2 u$).
In particular,  $u_{\ell-2} = \frac{\nu_\ell^2}{2} u_{\ell-1} > u_{\ell-1}$,
which implies that $s_{\ell-2}< 0$.
Now, it is not difficult to see from the fixed
point equation for $\gamma_{\nu_\ell}$
that for any $3\le k \le \ell$ we have 
$s_k= s_{k-1}/ f'(c_{k-1})$. It follows that
$s_{\ell-1}> 0$, and that $s_{2j}< 0$
for $4\le 2j \le \ell-2$ and $s_{2j+1}> 0$
for $3 \le 2j+1 \le \ell-3$.
In other words, $\gamma_{\nu_\ell}$ is increasing on
$(c_2, c_\ell)$ (with minimal value
$u_1$) and decreasing on $(c_\ell, c_1)$ (with minimal value $u_{\ell-1}=\nu_\ell u_1$).

Take a nonnegative $C^1$ function which is supported in $(c_{\ell-1}(\nu_\ell),1/2)$
for all $\ell$,
and note that $\int \varphi \gamma_{\sqrt 2} \, dx=u \int \varphi \, dx$.
Since $\int \varphi \gamma_{\nu_\ell} \, dx = u_2(\nu_\ell) \int \varphi \, dx$,
it suffices to show that there is a constant $K> 0$ so that
for all large enough $\ell$
$$
u_2(\nu_\ell) \le u -K \ell \nu_\ell^{-\ell} \, ,
$$
in order to prove the theorem.
Note that $|c_2(\nu_\ell)-c_2(\sqrt 2)|=O(\nu_\ell^{-\ell})$.
It follows that $\nu_\ell-\sqrt 2 =O(\nu_\ell^{-\ell})$
and that $u_2(\nu_\ell)-u_1(\nu_\ell)= (1-\nu_\ell^2/2) u_{\ell-1}/\nu_\ell
= O(\nu_\ell^{-\ell})$. Therefore, it is enough to
prove that
\begin{equation}\label{final}
u_1(\nu_\ell) \le u -K' \ell \nu_\ell^{-\ell} \, ,
\end{equation}
for some $K'> 0$ and all large enough $\ell$.

The rest of the proof is now similar to the argument in
Theorem~\ref{counter}.
Writing $(c_{m(j)}, c_{m'(j)})$ for the interval on
which $\gamma_{\nu_\ell}$ is constant equal to $u_j$, we have
\begin{align}\label{eqq}
1&= \sum_{j=1}^{\ell-1} u_j (c_{m' (j)}-c_{m(j)})\\
\nonumber &=
u_1(c_\ell-c_2) + \sum_{j=2}^{\ell/2} (u_j-u_1) (c_{m' (j)}-c_{m(j)})\\
\nonumber&\qquad\quad
+\nu_\ell u_1(c_1-c_\ell) +\sum_{j=\ell/2+1}^{\ell-2}
(u_j-\nu_\ell u_1)(c_{m' (j)}-c_{m(j)})\, .
\end{align}
If $3\le j \le \ell/2$ we have $u_j-u_1 \ge u_{j}-u_{j-1}> D \nu_\ell^{-\ell+2j}$
and $c_{m' (j)}-c_{m(j)}>D \nu_\ell^{-2j}$ for $D>0$ independent
of $\ell$ and $j$.
The case $j > \ell/2$ is similar.
It follows that there is $C' > 0$
so that the right-hand-side of (\ref{eqq}) is larger than
\begin{equation}\label{eqq'}
u_1(c_\ell-c_2) 
+\nu_\ell u_1(c_1-c_\ell)
+C' \ell \nu_\ell^{-\ell}=
u_1 d_\ell + C' \ell \nu_{\ell}^{-\ell}
\end{equation}
for all large enough $\ell$, where
\begin{equation}\label{dl}
d_\ell=\bigl ( c_\ell(\nu_\ell)-c_2(\nu_\ell) + 
\nu_\ell(c_1(\nu_\ell)-c_3(\nu_\ell) \bigr ) \, .
\end{equation}
We have thus proved that
$$
u_1 \le \frac{ 1- C' \ell \nu_\ell^{-\ell}}{d_\ell}\, .
$$
Combining the above bound with (\ref{noru}), (\ref{dd}), and
the easily proved fact that $|d-d_\ell|=O(\nu_\ell^{-\ell})$, we get
(\ref{final}).
\end{proof}


\begin{appendix}
\section {Uniform Lasota-Yorke estimates and spectral stability}
\label{uLY}

We  recall how to get uniform Lasota-Yorke estimates.
For $|t|< \epsilon$,
define $J_t:=(-\infty , f_t(c)]$ and $\chi_t:\real \to \{0,1,1/2\}$ by
$$
\chi(x)=\begin{cases}
0& x \notin J_t\\
1& x \in \mbox{int}\,  J_t\\
\frac{1}{2}& x=f_t(c)\, .
\end{cases}
$$
The two inverse branches of  $f_t$,
a priori defined on $[f_t(a),f_t(c)]$ and $[f_t(b),f_t(c)]$,
may be extended to $C^2$ maps $\psi_{t,+}: J_t \to (-\infty, c]$
and $\psi_{t,-}: J_t \to [c, \infty)$, with $\sup |\psi_{t,\sigma}'|<1$ for $\sigma=\pm$.
(in fact there is a $C^2$ extension of $\psi_{t,\pm}$ in a small neighbourhood
of $J_t$.) It is no restriction of generality
to assume that  $\psi_{t,+}(a_0)=a_0$ for all $t$.

Put
\begin{equation}\label{L1t}
\LL_{1,t} \varphi(x):=
\chi_t(x)
\psi'_{t,+}(x)\varphi(\psi_{t,+}(x))
+\chi_t(x)
|\psi'_{t,-}(x)| \varphi(\psi_{t,-}(x))\, .
\end{equation}
The first remark is that (see e.g.  \cite[Lemma 13]{Ke}) there is
$D\ge 1$ so that for any  $\varphi \in BV$, we have
\begin{equation}
\label{first'}
|\LL _{1,t} \varphi - \LL_1 \varphi |_1
\le D  |t| \|\varphi \|_{BV}\, ,  \forall |t| < \epsilon \, .
\end{equation}
Let $\lambda^{-1}<\inf_{x\ne 0} |f'(x)|$. Now, since $c$ is not periodic,
the proof of (3.26)  in \cite[p. 177]{Bb} yields
$D'\ge 1$  so that for all
small enough $|t|$
\begin{equation}\label{LYu}
\var \LL^m_{1,t} \varphi \le
\lambda^m \var \varphi + D' |\varphi|_1 \, , \forall m \ge 1\, .
\end{equation}

\section {Formal relation between $\Psi(1)$ and the  derivative of the SRB measure}
\label{susc}

We refer to \cite{Ruok} and \cite{BuL}, and references therein, for
uniformly hyperbolic instances where the SRB measure is smooth,
and where the susceptibility function  is related to its derivative.
In this appendix, we first recall (in our notation) Ruelle's {\it formal} argument 
\cite{Ru1ph} leading
to the consideration of $\Psi(1)$ as a candidate for the derivative.
We then give another (perhaps new) {\it formal} argument.
For simplicity, we consider only $X \equiv 1$.

The first step is rigorous:
By (\ref{LYu}) there are $C\ge 1$ 
and $\xi\in (\lambda, 1)$ that for all $|t|< \epsilon$
and all $k \ge 1$
\begin{equation}
|\int \varphi(x) \rho_t(x)\, dx - \int \varphi(f^k_t (x))\, dx|
\le C \xi^k \, .
\end{equation}
Now, since $\varphi$ is $C^1$, there is $s=s_k$ with
$|s|<|t|$ so that
\begin{equation}\label{taylor}
\varphi(f^k (x))- \varphi(f^k_t (x))
=t \sum_{n=0}^{k-1} \frac{d}{dy} \varphi(f^n_s(y)) |_{y=f^{k-n}_s(x)}\, .
\end{equation}
Next,
\begin{align}
\int \sum_{n=0}^{k-1} &\frac{d}{dy} \varphi(f^n_s(y)) |_{y=f^{k-n}_s(x)}\, dx\\
\nonumber
&= \int \sum_{n=0}^{k-1} \frac{d}{dy} \varphi(f^n_s(y)) \LL_{1,s}^{k-n}(\varphi_0)(y) \, dy\\
\nonumber &= \int  \sum_{n=0}^{k-1}  \varphi'(y) \LL_{0,s}^n(\LL^{k-n}_{1,s}(\varphi_0))(y) \, dy\, .
\end{align}
Letting (this is of course a formal step that is
not justified here) $t\to 0$ and $k \to \infty$ in the above formula, and using
that $\LL^m_1 (\varphi_0)\to \rho_0$ as $m \to \infty$, we would find
\begin{equation}
\frac{\int \varphi \rho_t \, dx-\int \varphi \rho_0 \, dx}{t} 
\sim \sum_{n=0}^\infty \int    \varphi'(y) \LL_0^n(\rho_0)(y) \, dy\quad
\mbox{ as } t \to 0\, ,
\end{equation}
as announced.

\medskip

Let us give now the second formal argument. Consider
$$
\int \varphi \frac{\rho-\rho_t}{t} \, dx
- \int \varphi'(\id -\LL_0)^{-1} \rho_0\, dx\, .
$$
Define for $x\in \real$ and $|t|< \epsilon$
\begin{equation}
R_t(x):=-1+\int_{-\infty}^x \rho_t(u) du\, .
\end{equation}

If $t$ is small, it is tempting (but of course illicit, since there
is no continuity of the resolvent on $BV$) to replace
$(\id -\LL_0)^{-1}$ by $(\id -\LL_{0,t})^{-1}$ where
\begin{equation}\label{L0t}
\LL_{0,t} \varphi(x):=
\chi_t(x) \varphi(\psi_{t,+}(x))-\chi_t(x)  \varphi(\psi_{t,-}(x))\, .
\end{equation}
We then get by
integration by parts
\begin{align}
\int \varphi' \biggl (\frac{R_0-R_t}{t} 
&- (\id -\LL_{0,t})^{-1} R'_0\biggr )\, dx\\
\nonumber &=\int \varphi' 
 (\id -\LL_{0,t})^{-1} \biggl (\frac{R_0-R_t}{t}-\LL_{0,t} (\frac{R_0-R_t}{t})-R'_0\biggr )\, dx\\
\nonumber &=
\int \varphi'(\id -\LL_{0,t})^{-1}\biggl ( \frac{R_0(\cdot)-R_0(\cdot-t)}{t}
-R'_0(\cdot)\biggr ) \, dx\\
\nonumber &\sim 0 \quad \mbox{ as } t \to 0\, ,
\end{align}
 where we used $\LL_{0,t} R_t=R_t$ and $\LL_{0,t} R_0(x)=\LL_0 R_0 (x-t)=R_0(x-t)$,
and where we ``pretend" again that $(\id -\LL_{0,t})^{-1}$ is continuous
on $BV$.

\section{A regularised susceptibility function}
\label{regs}

In this section we assume  (i)--(iv)
and, in addition, that $c$ is not preperiodic and $f(c)< b$,
$\min(f(a), f(b))> a> a_0$. 
For simplicity, we only consider the case $X \equiv 1$.

Define a power series
\begin{equation}
\rho_s(z)= \sum_{k=1}^\infty z^k s_k H_{c_k} \, .
\end{equation}
Clearly, (\ref{notper}) implies that $z\mapsto \rho_s(z)$ is holomorphic in a disc
of radius larger than $1$, with $\rho_s(1)=\rho_s$.
Put $\rho(z)=\rho_r+\rho_s(z)$.
Define next a regularised susceptibility function by
\begin{equation}
\tilde \Psi(z)=\sum_{n=0}^\infty \int  z^n \LL_0^{n} (\rho(z))(x) \varphi' (x)\, dx
\end{equation}

\begin{proposition}
Let $c$ be non preperiodic.
$\tilde \Psi(z)$ is holomorphic in the
open unit disc. If $\JJ(f)=0$ then
for every $\varphi \in C^1([a_0,b])$,  $\tilde \Psi(z)$
is holomorphic in a disc of radius larger than one,
and in addition,   $\tilde \Psi(1)=\Psi_1$.
\end{proposition}

\begin{proof}
The first claim is easily shown.
We have
for each $|z|< 1$
\begin{align}\label{singular}
\sum_{n=0}^\infty z^n \LL^n_0(\rho_s(z))  
&=\sum_{n=0}^\infty z^n \sum_{k=1}^\infty z^k s_{k} H_{c_{k+n}}\\
\nonumber &= \sum_{j=1}^\infty z^j H_{c_j} \sum_{k=1}^j  s_k \, .
\end{align}
Thus, since $\sup |H_{c_j}|\le 1$, (\ref{singular}) 
and (\ref{zerojump}) imply that $\tilde \Psi(z)$ is holomorphic in the disc of radius
$1/\xi$. The last claim  follows easily.
\end{proof}

Our second observation follows:

\begin{proposition}
Let $c$ be nonpreperiodic and $\JJ(f)\ne 0$.
For $\varphi \in C^1([a_0,b])$ with $\int \varphi \rho_0 dx =0$,
the limit of $\tilde \Psi(z)$ as $z\to 1$ in $[0,1)$ exists if and only
if 
\begin{equation}\label{notex}
\lim_{z \to 1, z \in [0,1)}\sum_{j=1}^\infty z^j \varphi(c_j) 
\end{equation}
exists.
\end{proposition}

\begin{proof}
Replacing $\rho_s$ by $\rho_s(z)$ in the proof of Proposition ~\ref{prop1},
it suffices to consider
$$
\sum_{n=1}^\infty \sum_{k=1} ^\infty s_k z^{k+n} H_{c_{k+n}}
=\sum_{j=1}^\infty z^j H_{c_j} \bigl (\sum_{k=1}^j s_k-\JJ(f) \bigr )
+ \JJ(f) \sum_{j=1}^\infty z^j H_{c_j} \, .  
$$
The first term in the right-hand-side of the above equality extends holomorphically in
the disc of radius $1/\xi$. Integrating by parts, this leaves
$$
\JJ(f) \sum_{j=1}^\infty z^j \varphi(c_j) \, ,
$$
as claimed.
\end{proof}
If $c_j$ is not recurrent it is easy to find examples of $\varphi$  so that
the limit (\ref{notex}) does not exist. 
This limit may never \cite{AA} exist.
\end{appendix}

\bibliographystyle{amsplain}

\begin{thebibliography}{10}

\bibitem{AA} A. Avila, personal communication, 2006.

\bibitem{Bb}
V. Baladi,
\textit{Positive transfer operators and decay of correlations,}
Advanced Series in Nonlinear Dynamics, Vol. 16,
Singapore, World Scientific,
2000.

\bibitem{BJR}
V. Baladi, Y. Jiang, and H.H. Rugh,
\textit{Dynamical determinants via dynamical conjugacies for postcritically finite polynomials,} J. Statist. Phys.  {\bf 108}  (2002) 973--993. 

\bibitem{BK}
V. Baladi and G. Keller,
\textit{Zeta functions and transfer operators for piecewise monotone transformations,}
Comm. Math. Phys. {\bf 127} (1990) 459--479.

\bibitem{BS} V. Baladi and D. Smania,
{\it work in progress} (2007).

\bibitem{BuL}
O. Butterley and C. Liverani,
\textit{Smooth Anosov flows: correlation spectra and stability,}
J. Modern Dynamics {\bf 1} (2007) 301--322.

\bibitem{Dolgo}
D. Dolgopyat, 
\textit{ On differentiability of SRB states for partially hyperbolic systems,}
Invent. Math. {\bf 155} (2004) 389--449. 

\bibitem{Er} S.V. Ershov, \textit{Is a perturbation theory for dynamical
chaos possible?,} Physics Letters A {\bf 177} (1993) 180--185.

\bibitem{Gal} G. Gallavotti,
\textit{Chaotic hypothesis: Onsager reciprocity
and fluctuation-dissipation theorem,}
J. Stat. Phys. {\bf 84} (1996) 899--926.

\bibitem{JL}
M. Jiang and R. de la Llave, 
\textit{ Linear response function for coupled hyperbolic attractors,}
  Comm. Math. Phys.  {\bf 261}  (2006) 379--404.

\bibitem{JR} Y. Jiang and D. Ruelle,
\textit{Analyticity of the susceptibility function for unimodal
markovian maps of the interval,}
 Nonlinearity  {\bf 18}  (2005) 2447--2453.

\bibitem{Ke} G. Keller,
\textit{Stochastic stability in some chaotic dynamical systems,}
Monatshefte Math. {\bf 94} (1982) 313--333.

\bibitem{Ke2}
G. Keller, 
\textit{On the rate of convergence to equilibrium in one-dimensional systems,}
Comm. Math. Phys. {\bf 96} (1984) 181--193.

\bibitem{Ke3} G. Keller,
\textit{An ergodic theoretic approach to mean field coupled maps,}
Progress in Probability, Vol. 46 (2000) 183--208.



\bibitem{KL}
G. Keller and C. Liverani,
\textit{Stability of the spectrum for transfer operators,} Annali  Scuola Normale Superiore di Pisa, 
{\bf 28} (1999) 141--152.



\bibitem{LY}
A. Lasota and J.A. Yorke, 
\textit{On the existence of invariant measures for piecewise monotonic transformations,}
Trans. Amer. Math. Soc. {\bf 186} (1973), 481--488 (1974).
 
\bibitem{LiY}
T.Y. Li and J.A. Yorke, 
\textit{Ergodic transformations from an interval into itself,}
Trans. Amer. Math. Soc. {\bf 235} (1978), 183--192. 

\bibitem{MM} M. Mazzolena, \textit{Dinamiche espansive unidimensionali:
dipendenza della misura invariante da un parametro,} Master's Thesis, Roma 2 (2007). 

\bibitem{RN} F. Riesz and B. Sz.-Nagy,
\textit{Functional analysis,} 1990 (reprint of the 1955 original),
Dover Books on Advanced Mathematics, Dover, New York.

\bibitem{Rub}
D. Ruelle, 
\textit{Dynamical zeta functions for piecewise monotone maps of the interval,}
CRM Monograph Series, 4. American Mathematical Society, Providence, RI, 1994.

\bibitem{Ru1} D. Ruelle,
\textit{Sharp zeta functions for smooth interval maps,}
In: International Conference on
Dynamical Systems (Montevideo 1995), Longman, 1996,
188--206.

\bibitem{Ru1ph} D. Ruelle,
\textit{General linear response formula in statistical mechanics, and the fluctuation-dissipation theorem far from equilibrium,}  Phys. Lett. A  {\bf 245}
(1998) 220--224.


\bibitem{Ruok}
D. Ruelle, \textit{Differentiation of SRB states,} 
Comm. Math. Phys. {\bf 187} (1997) 227--241. 

\bibitem{Ruok'}
D. Ruelle, \textit{Differentiation of SRB states: Corrections and complements,} 
Comm. Math. Phys. {\bf 234} (2003) 185--190. 

\bibitem{Rupr}
D. Ruelle, 
\textit{Application of hyperbolic dynamics to physics: some problems and conjectures,}  
Bull. Amer. Math. Soc.   {\bf 41}  (2004) 275--278.

\bibitem{Ru3}
D. Ruelle, \textit{Differentiating the a.c.i.m. of an interval map
with respect to $f$,}
Comm. Math. Phys. {\bf 258}  (2005) 445--453.


\bibitem{Ruwip} D. Ruelle, Work in progress
(personal communication, 2006)



\bibitem{RS}
M. Rychlik, and E. Sorets, 
\textit{Regularity and other properties of absolutely continuous invariant measures for the quadratic family,}
Comm. Math. Phys. {\bf 150} (1992) 217--236.

\bibitem{Sm} D. Smania, 
\textit{On the hyperbolicity of the period-doubling fixed point,}
Trans. Amer. Math. Soc. {\bf 358} (2006) 1827--1846. 

\bibitem{Ts}
M. Tsujii,
\textit{On continuity of Bowen--Ruelle--Sinai measures in families of one dimensional maps,}  Commun. Math. Phys. {\bf 177} (1996) 1--11


\bibitem{Yo}
L.-S. Young, 
\textit{What are SRB measures, and which dynamical systems have them?}
J. Statist. Phys. {\bf 108} (2002) 733--754
 
\end{thebibliography}

\end{document}